\documentclass[11pt,twoside]{article}
\usepackage{latexsym,amsmath}
\usepackage{indentfirst}
\usepackage{graphicx}

\numberwithin{equation}{section}
\newtheorem{Prop}{\bf Proposition}[section]
\newtheorem{Cor}{\bf Corollary}[section]

\newtheorem{Rem}{\bf Remark}[section]
\newtheorem{Ex}{\bf Example}[section]

\begin{document}
\def \b{\Box}
\begin{center}
{\Large {\bf A general method for stability controllability in the theory of fractional-order differential systems}}
\end{center}

\begin{center}
{\bf Gheorghe IVAN}\\
\end{center}

\setcounter{page}{1}
\pagestyle{myheadings}

{\small {\bf Abstract}. The main purpose of this paper is  to present a general method for the controllability of the stability  of a system of fractional-order differential equations around its equilibrium states. This method is applied to analyze and control the fractional stability of the fractional $ 2-$dimensional fractional Toda lattice with one linear control.}
{}

\section {Introduction}

The theory of fractional differential equations (i.e. fractional calculus) and its applications are based on  non-integer order of derivatives and integrals \cite{podl, kilb, bdst, sczc}. The fractional calculus has deep and natural connections with many fields of science and engineering.

In the last three decades, one increasing attention has been paid to the study of the dynamic behaviors (in particular, the chaotic behavior) of some classical differential systems, as well as some fractional-order differential systems. For example, the fractional models played an important role in  applied mathematics \cite{igim, ivan, migi}, mathematical physics \cite{nabu, puta, ivmi, kich, mihi}, applied physics \cite{ahma, sczc, mi22},  study of biological systems  \cite{ahme, pagi, lima}, chaos synchronization, secure communications \cite{zhli, bhda, igmp, muba} and so on.

The Lie groups, Lie algebroids and Leibniz algebroids have proven to be powerful tools for geometric formulation of the Hamiltonian mechanics \cite{demi, giop, gimo}. Also, they have been used  in the investigation of many fractional dynamical systems \cite{nabu, imod, migo}.

This paper is structured as follows. The Section 2 is devoted to the exposition  of the controllability method of fractional stability at an equilibrium state of a given fractional system. This method consists in associating a given fractional-order system with a new sistem of fractional-order differential equations, called the controlled fractional system around of an  equilibrium point. In Section 3 we investigate the fractional differential systems associated to
  2-dimensional Toda lattice with one linear control $ (3.4)~$ in terms of fractional Caputo derivatives. For this fractional-order model we investigate the existence and uniqueness of solution of initial value problem and asymptotic stability of its equilibrium states.
In Section $ 4,$ for the asymptotic stabilization  of fractional model  $ (3.4),$ we associate the controlled fractional-order system with  controls $~c_{1}, c_{2}~$ at the equilibrium point $~x_{e},~$ denoted by $ (4.4).~$  In Proposition $~(4.1)~$ are established sufficient conditions on parameters  $~k, c_{1}, c_{2}~$  to control the chaos in the fractional-order system $~(4.4).$\\[-0.5cm]

\section{Stability analysis and controllability of a fractional-order differential system}

We recall the Caputo definition of fractional derivatives, which is often used in concrete applications. Let $ f\in C^{\infty}(
\textbf{R}) $ and $ q \in \textbf{R}, q > 0. $ The $ q-$order Caputo differential operator
\cite{diet}, is described by $~D_{t}^{q}f(t) = I^{m - q}f^{(m)}(t), ~q > 0,~$ where $~f^{(m)}(t)$  represents the $
m-$order derivative of the function $ f,~m \in \textbf {N}^{\ast}$
is an integer such that $ m-1 \leq q \leq m $ and $ I^{q} $ is the
$ q-$order Riemann-Liouville integral operator, which
is expressed by $~I^{q}f(t)
=\displaystyle\frac{1}{\Gamma(q)}\int_{0}^{t}{(t-s)^{q
-1}}f(s)ds,~q > 0,$ where $~\Gamma $ is the Euler Gamma  function.
 If $ q =1$, then $ D_{t}^{q}f(t) = df/dt.~$

 In this paper we suppose that $ q \in (0,1].$
\markboth{Gh. Ivan}{ A general method for stability controllability in the theory ...}

 We consider the following system of fractional-order differential equations on ${\bf R}^{n}$:\\[-0.2cm]
\begin{equation}
D_{t}^{q}x^{i}(t) = f_{i}(x^{1}(t), x^{2}(t), \ldots, x^{n}(t)) ,~~ i=\overline{1,n},\label{(2.1)}
\end{equation}
where $ q\in (0,1), f_{i}\in C^{\infty}({\bf R}^{n}, {\bf R}),
~ D_{t}^{q} x^{i}(t)$ is the Caputo fractional derivative of order $ q $ for $ i=\overline{1,n}$ and $t\in [0,\tau)$ is the
time.

The fractional dynamical system $(2.1)$ can be written as follows:
\begin{equation}
D_{t}^{q}x(t) = f(x(t)),\label{(2.2)}
\end{equation}
where $~f(x(t)) = (f_{1}(x^{1}(t),\ldots, x^{n}(t)),
f_{2}(x^{1}(t),\ldots, x^{n}(t)), \ldots, f_{n}(x^{1}(t),\ldots,
x^{n}(t)))^{T} $ and $ D_{t}^{q} x(t)= ( D_{t}^{q}
x^{1}(t), \ldots, D_{t}^{q} x^{n}(t))^{T}.$

A point $ x_{e}=(x_{e}^{1}, x_{e}^{2},\ldots, x_{e}^{n})\in {\bf
R}^{n}$ is said to be {\it equilibrium state} of the fractional differential system
$(2.2)$, if $~D_{t}^{q}x^{i}(t) =0 $ for $ i=\overline{1,n}.~$. Its equilibrium states are determined by solving the set of equations: $~f_{i}(x^{1}(t),
x^{2}(t), \ldots, x^{n}(t)) = 0 ,~~ i=\overline{1,n}.$

The Jacobian matrix associated to $~(2.2) $ is $~J(x)=(\displaystyle\frac{\partial f_{i}}{\partial
x^{j}}),~~~i,j=\overline{1,n}.$

The stability of the fractional  system $(2.2)$ has been studied by Matignon
in \cite{mati}, where necessary and sufficient conditions have been established.
\begin{Prop} {\rm (\cite{mati})}
Let $ x_{e} $ be an equilibrium state of fractional differential system $(2.2)$ and
$ J(x_{e}) $ be the Jacobian matrix $J(x)$ evaluated at $ x_{e}$.

 $(i)~ x_{e}$ is locally asymptotically stable, if and only if  all eigenvalues
$ \lambda(J(x_{e})) $ of  $ J(x_{e}) $ satisfy:
\begin{equation}
| arg(\lambda (J(x_{e}))) | > \displaystyle\frac{q\pi}{2}.\label{(2.3)}
\end{equation}
 $(ii)~ x_{e} $ is locally stable, if and only if either it is asymptotically stable, or the
critical eigenvalues satisfying $~| arg(\lambda (J(x_{e}))) | = \displaystyle\frac{q \pi}{2}~$ have geometric
multiplicity one.\hfill$\Box$
\end{Prop}

Using the notation: $~~\tilde{q}:=\frac{2}{\pi}|arg(\lambda (J(x_{e})))|~$ and applying Proposition $(2.1)~$ one obtains the following corollary.
\begin{Cor}{\rm (\cite{danc})}
$(i)~$ The equilibrium state $~x_{e}~$  of the fractional model $ (2.2) $ is asymptotically stable if and only if the difference $~q-\tilde{q}~$ is strictly negative. More precisely, $~x_{e}~$ is asymptotically stable for all $~q\in (0, \tilde{q}).~$

$(ii)~$ If $~q-\tilde{q} > 0,~$ then  $~x_{e}~$ is unstable and the fractional model $ (2.2) $  may exhibit chaotic behavior. More precisely, $~x_{e}~$ is unstable $~(\forall) q\in (\tilde{q}, 1).~$ \hfill$\Box$
\end{Cor}
\begin{Cor} Let $~x_{e}~$ be an equilibrium state of the fractional model $~(2.2)~$ and $~\lambda_{i},~i=\overline{1,n}~$ the eigenvalues of $~J(x_{e}).$

$(i)~$ If one of the eigenvalues $~\lambda_{i},~i=\overline{1,n}~$ is equal to zero or it is positive, then $~x_{e}~$ is unstable  for all $~q\in (0,1).$

$(ii)~$ If $~\lambda_{i} < 0, $ for all $~i=\overline{1,n},~$ then  $~x_{e}~$   is asymptotically stable  $~(\forall)~q\in (0,1).$
\end{Cor}
{\it Proof.} $~(i)~$ We suppose $ \lambda_{1} \leq 0.~ $ Then $~\tilde{q}=0,~$  since $~arg(\lambda_{1} =0. $ We have  $~q-\tilde{q}=q > 0~$ and applying  Corollary  $~2.1(ii),~$ it follows that $~x_{e}~$ is unstable $~(\forall) q\in (0,1).$

$(ii)~$ Let $ \lambda_{i} < 0,~i=\overline{1,n}.~$ Then  $~|arg(\lambda (J(x_{e})))|= \pi~$ and $~ ~\tilde{q}=2.~$   We have  $~q-\tilde{q} < 0.~$ By Corollary  $~2.1(i),~$ it follows that $~x_{e}~$ is asymptotically stable $~(\forall) q\in (0,1).$ \hfill$\Box$

In the case when $ x_{e} $ is a unstable equilibrium state of the fractional differential system $ (2.2),$ we propose a simple method for to control
 the stability  of fractional model $ (2.2) $ around its equilibrium point $ x_{e}. $  This method consists in associating the fractional-order system $ (2.2) $ with a new sistem of fractional-order differential equations, determined by $ (2.2) $ and a control function $ u(x(t))\in C^{\infty}({\bf R}^{n}, {\bf R}).~$  To apply this method, the following four steps must be performed:\\
{\bf Step $~{\bf 1}.$}  We associate to $(2.2)$ a new fractional-order differential system defined by\\[-0.2cm]
\begin{equation}
D_{t}^{\alpha}x(t) = f(x(t))+ u(x(t)),\label{(2.4)}
\end{equation}
where  $~u(x(t)) = (u_{1}(x^{1}(t),\ldots, x^{n}(t)),\ldots,
x^{n}(t)))^{T}~$ is a control function.\\
{\bf Step $~{\bf 2}.$}  We choose the function control $~u(x(t)) $ such that  $~u(x_{e})= 0. $ \\
{\bf Step $~{\bf 3}.$} A good option for choosing the control function $~u(t)~$ that satisfies the relation $~u(x_{e})=0~$ is the following\\[-0.2cm]
   \begin{equation}
u_{i}(t) = c_{i} (x^{i}(t) - x_{e}^{i}),~~~ i=\overline{1,n}\label{(2.5)}
\end{equation}
where  $~c_{i}\in {\bf R} ~$ are control parameters.\\

With the control function $~u(t)~$  given by $~(2.6), $ the fractional system $~(2.4) $ becomes
\begin{equation}
D_{t}^{q}x^{i}(t) = f_{i}(x(t))+  c_{i} (x^{i}(t) - x_{e}^{i}),~~~ i=\overline{1,n},~~~~~ q\in (0,1)\label{(2.6)}
\end{equation}
where $~c_{i}\in {\bf R}, i=\overline{1,n}~.$

The  fractional model $~(2.5) $ is called the {\it controlled fractional system associated  to
fractional system $(2.2) $ at equilibrium point $~x_{e}.$}

{\bf Step $~{\bf 4}.$} For to analyse the fractional stability of the controlled fractional system $ (2.6) $ we apply the Matignon's test.

 This above method will be called the {\it controllability method of fractional stabilility at an equilibrium state of a given fractional system}.

It is easy to see that $ (2.2) $ and $ (2.6) $ have $~x_{e}~$ a common equilibrium state.
\begin{Rem}{\rm $(i)~$ This method was applied to analyse and control the fractional stability type of steady states for the fractional differential equations $~3D~$ Maxwell-Bloch type in \cite{ghiv}.

$(ii)~$ Other fractional modeling of classical dynamical systems in which this method is applied to control their fractional stability have been discussed in \cite{ivmi, miha, migi}.} \hfill$\Box$\\[-0.4cm]
\end{Rem}

If one selects the appropriate parameters $c_{i}, i=\overline{1,n}$  which then make the eigenvalues of the
linearized equation of $(2.6)$ satisfy one of the conditions from Proposition 2.1, then the trajectories of $ (2.6) $ asymptotically
approaches the unstable equilibrium state $x_{e}$ in the sense that $\lim_{t\rightarrow \infty} \|x(t)-x_{e}\|= 0$, where
$\|\cdot\|$ is the Euclidean norm. In the case when the equilibrium state $~x_{e}~$ is unstable, then fractional model $ (2.6) $  may exhibit chaotic behavior.\\[-0.5cm]

\section{ Stability analysis of the $ 2-$dimensional fractional-order Toda lattice with one linear control}

The Toda-type systems \cite{dami} are described by the following equations on $ {\bf R}^{2n-1}: $\\[-0.2cm]
\begin{equation}
 \dot{x}^{i}(t)  =  x^{i}(t)( y^{i+1}(t)- y^{i}(t)), ~~~~~
 \dot{y}^{j}(t)  =  2[(x^{j})^{2}(t) - (x^{j-1})^{2}(t)], \label{(3.1)}
\end{equation}
where $~x^{0}(t) = x^{n}(t)=0,~x^{i}, i=\overline{1,n-1}, ~y^{j}, j=\overline{1,n}~$ are state variables, $~\dot{x}^{i}(t)= dx^{i}(t)/dt,~\dot{y}^{j}(t)= dy^{j}(t)/dt~$ and $ t $ is the time. The system $ (3.1)$  is called the {\it n-dimensional Toda lattice}.

The {\it n- dimensional fractional-order  Toda lattice } associated to dynamics $~(3.1)~$ is defined by the following
set of fractional differential equations:\\[-0.2cm]
\begin{equation}
\left\{ \begin{array} {lll}
 D_{t}^{q}{x}^{i}(t) & = & x^{i}(t)( y^{i+1}(t)- y^{i}(t)), ~~~~~~~~ i=\overline{1,n-1}  \\[0.1cm]
 D_{t}^{q}{y}^{j}(t) & = & 2[(x^{j})^{2}(t) - (x^{j-1})^{2}(t)],~~~~~j=\overline{1,n}~~~  q \in (0,1),\\[0.1cm]
  x^{0}(t) = 0, & & x^{n}(t)=0. \label{(3.2)}
  \end{array}\right.
\end{equation}
In this section we investigate the n- dimensional fractional-order Toda lattice for $~n=2 $ with one linear control about $~Oy^{2}-$axis. This fractional model is described  by:\\[-0.2cm]
\begin{equation}
 D_{t}^{q} x^{1} =  x^{1}(-y^{1} +  y^{2}),~~~ D_{t}^{q} y^{1}  = 2 (x^{1})^{2},~~~
 D_{t}^{q} y^{2}  = -2 (x^{1})^{2}- k y^{2}, \label{(3.3)}
\end{equation}
where $~k\in {\bf R}^{*}~$ is a control parameter.

Using the transformations $ x^{1} = x^{1}, ~ y^{1} = x^{2},~ y^{2}
= x^{3}, $ the system $(3.3)$ becomes:\\[-0.2cm]
\begin{equation}
\left\{ \begin{array} {lll}
 D_{t}^{q}{x}^{1}(t) & = & x^{1}(t)(- x^{2}(t) + x^{3}(t)), \\[0.1cm]
 D_{t}^{q}{x}^{2}(t) & = & 2 (x^{1}(t))^{2},          ~~~~~~~~~~~~~~~~~  q \in (0,1),\\[0.1cm]
 D_{t}^{q}{x}^{3}(t) & = & - 2 (x^{1}(t))^{2}- k x^{3}(t). \label{(3.4)}
  \end{array}\right.
\end{equation}
The initial value problem of the fractional system
$(3.4)$ can be represented in the following matrix form:\\[-0.2cm]
\begin{equation}
D_{t}^{\alpha}x(t)  =  x^{1}(t) A x(t) + x^{3}(t) B x(t),~~~~~~~~ x(0) =
x_{0},\label{(3.5)}
\end{equation}
where $0 < q < 1,~ x(t)= ( x^{1}(t),
 x^{2}(t), x^{3}(t))^{T}, ~t\in(0,\tau)$ and\\[-0.1cm]
\[
A = \left ( \begin{array}{ccc}
0 & -1 & 1 \\
2 & 0 & 0 \\
-2 & 0 & 0 \\
\end{array}\right ),~~~ B = \left ( \begin{array}{ccc}
0 & 0 & 0 \\
0 & 0 & 0\\
0 & 0 & -k\\
\end{array}\right ).
\]
\begin{Prop}
The initial value problem of the 2-dimensional fractional-order Toda lattice with one control
$(3.4)$ has a unique solution.
\end{Prop}
{\it Proof.} Let $ f(x(t))= x^{1}(t) A x(t) + x^{3}(t) B x(t). $ It is
obviously continuous and bounded on $ D =\{ x \in {\bf R}^{3} |~
x^{i}\in [x_{0}^{i} - \delta, x_{0}^{i} + \delta]\}, i=\overline{1,3} $ for any
$\delta>0. $\\
 We have $~f(x(t)) - f(y(t)) =  x^{1}(t) A x(t) - y^{1}(t) A y(t) + x^{3}(t) B x(t)-  y^{3}(t) B y(t)=g(t)+h(t), $
where $~g(t)= x^{1}(t) A x(t) - y^{1}(t) A y(t)~$ and $~h(t)= x^{3}(t) B x(t) - y^{3}(t) B y(t).~$ Then\\[0.1cm]
 $(a)~~|f(x(t)) - f(y(t))|\leq |g(t)| + |h(t)|. $\\
Using reasoning analogous to that in the proof of the Proposition 2.1 in \cite{ivmi}, we can show that:\\[0.1cm]
$(b)~~~|g(t)| \leq (\|A\|+ |y^{1}(t)|)\cdot |x(t)- y(t)|~~$   and  $~~~|h(t)| \leq (\|B\| + |y^{3}(t)|)\cdot |x(t)- y(t)|.$\\[0.1cm]
According to $(b)~$ the relation $(a)$ becomes\\[0.1cm]
$(c)~~~|f(x(t)) - f(y(t))|\leq  (\|A\| + \|B\| + |y^{1}(t)| + |y^{3}(t)|)\cdot |x(t)-y(t)|.$\\[0.1cm]
Replacing  $\|A\|= \sqrt{10},~ \|B\| = |k|~$ and using the inequalities $~|y^{i}(t)|\leq |x_{0}| + \delta, ~i=1,3~$ from the
relation $ (c), $ we deduce that\\[0.1cm]
$(d)~~~|f(x(t)) - f(y(t))|\leq  L\cdot |x(t)-y(t)|,~~~~~
\hbox{where}~ L =\sqrt{10} + |k| + 2(|x_{0}| + \delta) > 0.$\\[0.1cm]
 The inequality $(d)$ shows that $ f(x(t))$
satisfies a Lipschitz condition. Based on the results of Theorems
$1$ and $2$ in \cite{difo}, we can conclude that the initial value
problem of the system $(3.3)$ has a unique solution. \hfill$\Box$

For the fractional system $ (3.4) $ we introduce the following notations:\\[-0.2cm]
 \begin{equation}
f_{1}(x) =  - x^{1}x^{2} + x^{1}x^{3},~~~ f_{2}(x) = 2 (x^{1})^{2},~~~ f_{3}(x)  = - 2 (x^{1})^{2} - k x^{3}.\label{(3.5)}
\end{equation}
\begin{Prop}
{\it The equilibrium states of the 2-dimensional fractional-order Toda lattice
$(3.4)$ are given as the following family}:\\[-0.3cm]
\[
E:=\{ e_{m} =(0, m, 0)\in {\bf R}^{3} |~ m \in {\bf R}\}.
\]
\end{Prop}
{\it Proof.} The equilibrium states are solutions of the equations
$~f_{i}(x)=0, i=\overline{1,3}$ where $~f_{i},~i=\overline{1,3}$
are given by (3.6).\hfill$\Box$

Let us we present the study of asymptotic stability of equilibrium
states for the fractional system $(3.4)$. Finally, we will discuss
how to stabilize the unstable equilibrium states of the system
$(3.4)$ via fractional order derivative. For this study we apply
the Matignon's test.

The Jacobian matrix associated to system $(3.4)$ is:
\[
J(x,k) = \left ( \begin{array}{ccc}
-x^{2}+x^{3} & -x^{1}    & x^{1} \\
  4 x^{1}   & 0 & 0\\
 -4 x^{1}   &  0  & -k \\
\end{array}\right ).\\[-0.1cm]
\]
\begin{Prop}
The equilibrium states $ e_{m}\in E~$
are unstable $ (\forall) q \in (0,1).$
\end{Prop}
{\it Proof.} The characteristic polynomial of the matrix $~
J(e_{m},k) =\left (\begin{array}{ccc}
  -m & 0  & 0\\[0.1cm]
  0 & 0 & 0 \\
  0 & 0  & -k \\
\end{array}\right ) $
is\\
 $~ p_{J(e_{m},k)}(\lambda) = \det ( J(e_{m},k) -
\lambda I) = - \lambda (\lambda + m)(\lambda + k).~$
 The equation $ ~p_{J(e_{m},k)}(\lambda) = 0 $  has the root $ \lambda_{1} = 0.~$ By Corollary  2.2(i), follows that $ e_{m}, m\in {\bf R}~$ are unstable for
all $ q\in (0,1).$ \hfill$\Box$\\[-0.5cm]

\section{ Controllability of chaotic behavior of the fractional  model $~(3.4)~$}

For the controllabilty of chaotic behavior of the $~2-$dimensional Toda lattice with one linear control, we apply the controllability method of fractional stabilility of the fractional  model $~(3.4)~$ at an equilibrium point.

Let  $ x_{e}$ be an unstable equilibrium state. We associate to $(3.4)$ a new fractional-order system with (external) controls and given by:\\[-0.2cm]
\begin{equation}
\left\{ \begin{array} {lll}
 D_{t}^{q}{x}^{1}(t) & = & x^{1}(t)(- x^{2}(t) + x^{3}(t)) + u_{1}(t), \\[0.1cm]
 D_{t}^{q}{x}^{2}(t) & = & 2 (x^{1}(t))^{2}  + u_{2}(t),          ~~~~~~~~~~~~~~~~~  q \in (0,1),\\[0.1cm]
 D_{t}^{q}{x}^{3}(t) & = & - 2 (x^{1}(t))^{2}- k x^{3}(t)  + u_{3}(t), \label{(4.1)}
  \end{array}\right.
\end{equation}
where $ u_{i}(t), i=\overline{1,3}$ are control functions.

In this section we take the control functions $ u_{i}(t), i=\overline{1,3}, $ given by:\\[-0.2cm]
 \begin{equation}
u_{1}(t) = c_{1}(x^{1}(t)- x_{e}^{1}),~~~ u_{2}(t) = c_{2} (x^{2}(t)- x_{e}^{2}),~~~ u_{3}(t)  = 0, ~~~ c_{1}, c_{2} \in {\bf R}^{*}.\label{(4.2)}
\end{equation}
With the control functions $(4.2), $ the system $(4.1) $ becomes:\\[-0.2cm]
\begin{equation}
\left\{ \begin{array} {lll}
 D_{t}^{q}{x}^{1}(t) & = & x^{1}(t)(- x^{2}(t) + x^{3}(t)) + c_{1}(x^{1}(t) - x_{e}^{1}) , \\[0.1cm]
 D_{t}^{q}{x}^{2}(t) & = & 2 (x^{1}(t))^{2}  + c_{2} (x^{2}(t)-x_{e}^{2}),          ~~~~~~~~~~~~~~~~~  q \in (0,1),\\[0.1cm]
 D_{t}^{q}{x}^{3}(t) & = & - 2 (x^{1}(t))^{2}- k x^{3}(t), \label{(4.3)}
 \end{array}\right.
\end{equation}
where  $~k, c_{1}, c_{2}\in{\bf R}^{\ast}~$ are control parameters.

The fractional system $ (4.3) $ is called the {\it controlled fractional-order system associated  to
$~(3.4) $ at  $~x_{e}.$}

The controlled fractional-order system  associated to $ (3.4)~$ at  $~x_{e} = e_{m},$ is written:\\[-0.2cm]
\begin{equation}
\left\{ \begin{array} {lll}
 D_{t}^{q}{x}^{1}(t) & = & x^{1}(t)(- x^{2}(t) + x^{3}(t)) + c_{1}x^{1}(t) , \\[0.1cm]
 D_{t}^{q}{x}^{2}(t) & = & 2 (x^{1}(t))^{2}  + c_{2} (x^{2}(t)- m),          ~~~~~~~~~~~~~~~~~  q \in (0,1),\\[0.1cm]
 D_{t}^{q}{x}^{3}(t) & = & - 2 (x^{1}(t))^{2}- k x^{3}(t), \label{(4.4)}
 \end{array}\right.
\end{equation}
where  $~k, c_{1}, c_{2}\in{\bf R}^{\ast}~$ are control parameters and $~m\in {\bf R}.$

The Jacobian matrix of the fractional model $(4.4)$  is\\[-0.2cm]
\[
J(x, k, c_{1}, c_{2}) = \left (\begin{array}{ccc}
-x^{2}+x^{3}+c_{1} & -x^{1}    & x^{1} \\
  4 x^{1}   & c_{2} & 0\\
 -4 x^{1}   &  0  & -k \\
\end{array}\right ).\\[-0.1cm]
\]
\begin{Prop}
Let  be the fractional system  $(4.4)~$ and $~e_{m}= (0,m,0)\in E.$\\
$~~~~~~~$ {\bf 1.} $~k > 0.$\\
$(i)~$ If $~ c_{2}< 0,~$ then $ e_{m}~$ is asymptotically stable for all $~m\in (c_{1}, \infty )~$ and $ ~q\in (0,1).$\\
$(ii)~$ If $~c_{2} < 0,~$ then $~ e_{m}~$ is unstable for all $~m\in (-\infty, c_{1}]~$ and $~q\in (0,1).$\\
$(iii)~$ If $~c_{2} > 0,~$ then $~ e_{m}~$ is unstable $~(\forall) m\in {\bf R}~$ and $~q\in (0,1).$\\
$~~~~~~~$ {\bf 2.} $~k < 0.~$ If  $~c_{1}, c_{2}\in {\bf R}^{\ast},~$ then $~ e_{m} $ is unstable for all $~m\in {\bf R}~$ and $~q\in (0,1).$
\end{Prop}
{\it Proof.} The characteristic polynomial of matrix $~
 J(e_{m},k, c_{1}, c_{2}) =\left (\begin{array}{ccc}
  -m + c_{1} & 0  & 0\\[0.1cm]
  0 & c_{2} & 0 \\
  0 & 0  & -k \\
\end{array}\right )~$\\\
 is $~ p_{J(e_{m},k, c_{1}, c_{2})}(\lambda) = \det ( J(e_{m},k, c_{1}, c_{2}) -
\lambda I) = - (\lambda +m-c_{1}) (\lambda - c_{2})(\lambda + k).~$
The roots of equation $~ p_{J(e_{m},k, c_{1}, c_{2})}(\lambda) = 0~$ are $\lambda_{1}= c_{1}-m,~\lambda_{2}=c_{2},~\lambda_{3}=-k.$\\
{\bf 1.} {\bf Case} $~k >0~$ and $~q\in (0,1).~$  Then $~\lambda_{3} <0.~$\\
$(i)~$ We have $~\lambda_{1}<0~$ and $~\lambda_{2}<0~$  if and only if $~c_{2}<0~$  and  $~m\in (c_{1}, \infty).~$ It follows that $~\lambda_{i} <0, i=\overline{1,3}~$ and according to Corollary 2.2(ii), the equilibrium state $~ e_{m}~$ is asymptotically stable.\\
$(ii)~$ We suppose $~c_{2} < 0.~ $  Then $~\lambda_{2} <0.~$  We have $~\lambda_{1}>0~$ if and only if $~m\in (-\infty, c_{1}].~$  In this case, $~J(e_{m},k, c_{1}, c_{2})~$ has a positive eigenvalue and by Corollary  2.2(i), it follows that  $ e_{m} $ is unstable.\\
$(iii)~$ We suppose $~c_{2} > 0.~$ In this case, one from the eigenvalues $~\lambda_{1}~$  and $~\lambda_{2}~$ is positive. Since $~J(e_{m},k, c_{1}, c_{2})~$ has at least a positive eigenvalue, it follows that  $ e_{m} $ is unstable for all $~m\in {\bf R}.$ \\
{\bf 2.} {\bf Case} $~k <0~$ and $~q\in (0,1).~$ Then $~\lambda_{3}>0.~$ Since $~J(e_{m},k, c_{1}, c_{2})~$ has at least a positive eigenvalue, it follows that  $ e_{m} $ is unstable $~(\forall) m\in {\bf R}.~$ Hence, the assertions  $~(1)~$ and $~(2)~$ hold. \hfill$\Box$
\begin{Ex}
{\rm  Let be the 2-dimensional fractional-order Toda lattice with controls $k, c_{1}, c_{2}~$ described by $(4.4).~$\\
$(1)~$ We select $~ k=0.4, c_{1}  = - 0.02~$ and  $~c_{2}= -0.3.~$ According to Corollary $~2.2(ii),~$ it follows that  $~ e_{0}=(0,0,0)~$ is
asymptotically stable for all $~ q\in (0, 1).$\\
$(2)~$ We consider $~ k= 1, c_{1}  = 0.02~$ and  $~c_{2}= -0.3.~$ Applying Corollary $~2.2(i),~$ it follows that $~ e_{0}=(0,0,0)~$ is unstable for all $~ q\in (0, 1).~$ In other words, the fractional model $~(3.3)~$ with $~ k= 1, c_{1}  = 0.02~$ and  $~c_{2}= -0.3,~$ behaves chaotically around the equilibrium point $~ e_{0}.$\\
$(3)~$ We consider $~ k=-0.25, c_{1}  =- 0.02~$ and  $~c_{2}= -0.3.~$ Applying Corollary $~2.2(i),~$ it follows that $~ e_{0}=(0,0,0)~$ is unstable for all $~ q\in (0, 1).$}
\end{Ex}
\begin{Rem}{\rm  Toda-type dynamic systems have been studied from various research directions by many authors. More specifically, from the point of view of Poisson geometry, Toda lattices were discussed in the papers \cite{dami, pudo, ptud, puwa}, and as fractional-order differential systems they were investigated in \cite{imgo, igom}.} \hfill$\Box$
\end{Rem}
{\bf Conclusions.} This paper presents the 2-dimensional fractional-order Toda lattice with one liniar control, denoted by $ (3.4).$  The asymptotic stability of equilibrium states of $~(3.4) $ was investigated. Sufficient conditions on the parameters $~k, c_{1}, c_{2}~$ in the fractional model $(4.4)$ so that its equilibrium states are asymptotically stable have proved.  By choosing the right parameter $~k,c_{1}, c_{2}~$ in  $~(4.4),~$ this work offers a series of chaotic fractional differential systems.

 {\bf Acknowledgments.} The author wishes to tank the referees for their useful comments and suggestions.\\[-0.5cm]

Author's adress\\[-0.5cm]

Gheorghe Ivan\\[0.1cm]
West University of Timi\c soara. Seminarul de Geometrie \c si Topologie.\\
Department of Mathematics. Timi\c soara, Romania.\\
E-mail: gheorghe.ivan@e-uvt.ro\\

\begin{thebibliography}{99}
\smallskip
\bibitem{ahma} W.M. Ahmad, J.C. Sprott, {\it Chaos in fractional order autonomous nonlinear systems}. Chaos, Solitons and Fractals, {\bf 16}(2003), 339–351.
\bibitem{ahme} E. Ahmed , A.M.A. El-Sayed, H.A.A. El-Saka, {\it Equilibrium points, stability and numerical solutions of fractional order predator–prey and rabies models}. J. Math. Anal. Appl., {\bf 325}(2007), no. 1, 542–553. DOI:10.1016/j.jmaa.2006.01.087.
\bibitem{bdst} D. B\u aleanu, K. Diethelm, E. Scalas, J.J. Trujillo, {\it Fractional Calculus Models and Numerical Methods}. Series on Complexity, Nonlinearity and Chaos. World Scientific, 2012.
\bibitem{bhda} S. Bhalekar, V. Daftardar-Gejji, {\it Synchronization of different fractional order chaotic systems using active control}. Commun. Nonlinear Sci. and Numer. Simulat., {\bf 15}(2010), no. 11,  3536–3546. DOI:10.1016/j.cnsns.2009.12016.
\bibitem{dami} P.A. Damianou, {\it Multiple Hamiltonian structures for Toda-type systems }, J. Math. Phys., {\bf 35} (1994), 5511- 5541.
\bibitem{danc} M.F. Danca, {\it Hidden chaotic attractors in fractional-order systems}. Nonlinear Dynamics, {\bf 89} (2017), no. 1, 577-586. DOI:10.1007/s11071-017-3472-7.
\bibitem{demi}  M. Degeratu, M. Ivan, {\it Linear connections on Lie algebroids}, Proceedings of the 5th Conference of Balkan Society of Geometers, Aug. 29- Sept. 2, 2005. Mangalia Romania. Geometry Balkan Press, 2006, 44-53.
\bibitem{diet} K. Diethelm, {\it  The Analysis of Fractional Differential Equations: An Application-Oriented Exposition  Using Differential Operators of Caputo Type}. Springer, 2010.
 \bibitem{difo} K. Diethelm, N. J. Ford, {\it Analysis of fractional differential equations}. J. Math. Anal. Appl., {\bf 265}(2002), 229–248.
\bibitem{nabu} R.A. El-Nabulsi, {\it A fractional action-like variational approach of some classical quantum and geometrical dynamics}. Int. J. Appl. Math.,{\bf
    17}(2005), 299-317.
\bibitem{ivan} G. Ivan, {\it Geometrical and dynamical properties of general Euler top system}. Indian J. Pure Appl. Math., {\bf 44}(2013), no. 1, 77-93.
    DOI:10.1007/s13226-013-0004-0.
\bibitem{ghiv} G. Ivan, {\it On fractional differential equations of 3D Maxwell-Bloch type}. Int. J. Geom. Met. Mod. Phys., {\bf 11}(2014), no. 4, 1450028 (12 pages).  DOI:10.1142/S02198878.14500285
\bibitem{igim} G. Ivan, M. Ivan, {\it General Euler top system and its Lax representation}.  Int. J. Geom. Met. Mod. Phys., {\bf 8}(2011), no. 5, 937-944.
    DOI:10.1142/S0219887811003543.
\bibitem{giop} G. Ivan, D. Opri\c s, {\it Dynamical systems on Leibniz algebroids}, Diff. Geometry-Dynamical systems, {\bf 8}(2006). Geometry Balkan Press, 127-137.
\bibitem{gimo}  G. Ivan, M. Ivan, D. Opri\c s, {\it Fractional Euler-Lagrange and fractional Wong equations for Lie algebroids}, Proceed. of The 4-th Int. Colloq. "Math. in Eng. and Numerical Phys." October 6-8, 2006, Bucharest, Romania, 73-80.
\bibitem{imod} G. Ivan, M. Ivan, D. Opri\c s, {\it Fractional dynamical systems on fractional Leibniz algebroids}, Analele \c Stiin\c tifice ale Universit\u a\c tii "Al. I. Cuza" din Ia\c si (S.N.), Matematic\u a, {\bf 53} (2007), Supl., 222-234.
\bibitem{igmp} G. Ivan, M. Ivan, C. Pop, {\it Numerical integration and synchronization for the 3-dimensional metriplectic Volterra system}. Math. Probl. Eng., {\bf 2011}, Article ID 723629 (11 pages). DOI:10.1155/2011/723629.
\bibitem{igom} G. Ivan, D. Opri\c s,  M. Ivan, {\it Hybrid fractional differential systems associated to Toda lattice, Numerical simulation}. Proc. of the $12^{th}$ Symposium of Math. and Its Appl. "Politehnica" University of Timi\c soara, November 5-7, 2009, 227-234.
\bibitem{ivmi} M. Ivan, {\it Control chaos in the fractional Lorenz-Hamilton system}. Fractional Differ. Calc., {\bf 6}(2016), no. 1, 111-119. DOI:10.7153/fdc-06-07.
\bibitem{miha} M. Ivan, {\it Stability analysis and control chaos for fractional $5D$ Maxwell-Bloch model}. Preprint arXiv math.,2018. DOI:10.48550/arxiv.1802.07706.
\bibitem{mi22} M. Ivan,  {\it A fractional model for the single Stokes pulse from the nonlinear optics}, Journal of Applied Mathematics and Physics, {\bf 10}(2022), 2856-2875. DOI: 10.4236/jamp.2022.1010191.
\bibitem{mihi} M. Ivan, {\it Dynamics analysis of the fractional-order Lagrange system}.  Int. J. Modern Eng. Research, {\bf 12}(2022), no.8, 23-31.
\bibitem{migi} M. Ivan, G. Ivan, {\it On the fractional Euler top system with two parameters}.  Int. J. Modern Eng. Research, {\bf 8}(2018), no. 4, 10-22.
\bibitem{imgo} M. Ivan, G. Ivan,  D. Opri\c s, {\it Fractional differential systems associated to Toda lattice}, BSG Proceed. 16. The Int. Conf. of Diff. Geom. and Dynamical Systems (DGDS-2008) and The V-th Int. Colloq. of Math. in Eng. and Numerical Phys. (MENP-5), Aug. 29 – Sept. 2, 2008, Mangalia, Romania, pp. 80-90.
\bibitem{migo} M. Ivan, G. Ivan, D. Opri\c s, {\it Fractional equations of the rigid body on the pseudo-orthogonal group $ SO(2,1) $}. Int. J. Geom. Met. Mod. Phys., {\bf 6}(2009), no. 7, 1181-1192.  DOI:10.1142/S0219887809004168.
\bibitem{lima} J. Llibre, Y.P. Martinez, {\it Dynamics of a family of Lotka-Volterra Systems in $~{\bf R}^{\ast}. $} Nonlinear Anal., {\bf 199}(2020),111915.
\bibitem{kilb} A.A. Kilbas, H.M. Srivastava, J.J. Trujillo, {\it Theory and Applications of Fractional Differential Equations}. North-Holland, Math. Studies, {\bf 204}, 2006.
\bibitem{kich} P. Kumar, S.K. Chaudhary, {\it Analysis of fractional order control system with performance and stability}. Int. J. Eng. Sci. Tech. {\bf 9} (2017), no. 5, 408–416.
\bibitem{mati} D. Matignon, {\it Stability results for fractional differential equations with applications to control processing}. In: Proc. of the Comput. Eng. in Systems and Applications, IMACS, IEEE-SMC, Lille, France, July 1996, {\bf 2}(1996), 963-968.
 \bibitem{muba} P. Muthukumar, P. Balasubramaniam, {\it Feedback synchronization of the fractional order reverse butterfly-shaped chaotic system and its application to digital cryptography}. Nonlinear Dyn., {\bf 74}(2013), no. 4, 1169–1181. DOI:10.1007/s11071-013-1032-3.
\bibitem{podl} I. Podlubny, {\it Fractional Differential Equations}. Academic Press, New York, 1999.
\bibitem{pagi} C. Pop, A. Aron, C. Galea, M. Ciobanu, M. Ivan, {\it Some geometric aspects in theory of Lotka-Volterra system}. Proc. of the 11-th WSEAS
Int. Conf. on Sust. in Sci. Eng., Timi\c soara, Romania, May (2009), 91-97.
\bibitem{puta} M. Puta, C. Pop, C. D\u an\u aiasa, C. Hedrea, {\it Some geometric aspects in the theory of Lagrange system}, Tensor, N.S. {\bf 69}(2008), 83-87.
\bibitem{pudo} M. Puta, R. Tudoran, {\it Controllability, stability and the n-dimensional Toda lattice}, Bull. Sci. math. 126 (2002) 241–247.
\bibitem{ptud} M. Puta, Ra. Tudoran, R. Tudoran, {\it Poisson manifolds and Bermejo-Fairen construction of Casimirs}, Tensor, N.S. {\bf 66}(2005), 59-70.
\bibitem{puwa} M. Puta, D. Wainberg, {\it The stability of the equilibrium states for some mechanical systems}, Studia Univ. "Babe\c s - Bolyai", Mathematica, {\bf 54} (2009), no.1, 119-126.
\bibitem{sczc} H.G. Sun, A. Chang, Y. Zhang, W. Chen, {\it A review on variable-order fractional differential equations: mathematical foundations, physical models, numerical models and applications}. Fract. Calc. Appl. Anal., {\bf 22}(2019), no. 1, 27–59. DOI: 10.1515/fca-2019-0003.
\bibitem{zhli} T. Zhou, C. Li, {\it Synchronization in fractional-order differential systems}. Physica D, Nonlinear Phenomena, {\bf 212} (2005), 111-125. DOI:10.1016/j.physd.2005.09.012.
\end{thebibliography}
\end{document}